\documentclass{amsart}
\usepackage{fullpage}

\newcommand{\R}{{\mathbb R}}

\newcommand{\Mi}[1]{M^{[#1]}}

\newtheorem{theorem}{Theorem}
\newtheorem{lemma}[theorem]{Lemma}
\newtheorem{prop}[theorem]{Proposition}

\newtheorem*{theorem*}{Theorem}

\theoremstyle{remark}

\title{Partial monoids and Dold-Thom functors}

\author{Jacob Mostovoy}

\address{Departamento de Matem\'aticas, CINVESTAV \\
 Apartado Postal 14-740
C.P.\ 07000, M\'exico, D.F., M\'exico\\
E-mail: {\tt jacob@math.cinvestav.mx}}

\begin{document}

\begin{abstract}
Dold-Thom functors are generalizations of infinite symmetric products, where integer multiplicities of points are replaced by composable elements of a partial abelian monoid. It is well-known that for any connective homology theory, the machinery of $\Gamma$-spaces produces the corresponding linear Dold-Thom functor. In this note we construct such functors directly from spectra by exhibiting a partial monoid corresponding to a spectrum.
\end{abstract}

\maketitle


\section{Introduction.}

Let $SP^{\infty}X$ be the infinite symmetric product of a pointed connected cell complex $X$. Then, according to the Dold-Thom Theorem \cite{DT}, the homotopy groups of $SP^{\infty}X$ coincide, as a functor, with the reduced singular homology of $X$. Although there is no computational advantage in this definition of singular homology, it is important since it can be extended to the algebro-geometric context; in particular, the motivic cohomology of \cite{V} is defined in this way.

The construction of the infinite symmetric product can be generalized so as to produce an arbitrary connective homology theory. Such generalized symmetric products were defined by G.\ Segal in \cite{Segal:Topology}: essentially, these are labelled configuration spaces, with labels in a $\Gamma$-space (see also \cite{And, BF, Ku}). If the $\Gamma$-space of labels is injective
(see \cite{SV}) it gives rise to a partial abelian monoid; it has been proved that in \cite{SV} that for each connective homology theory there exists an injective $\Gamma$-space. In this case Segal's generalized symmetric product can be thought of as a space of configurations of points labelled by composable elements of a partial monoid. An explicit example discussed in \cite{Segal:K-theory} is the space of configurations of points labelled by orthogonal vector spaces: it produces connective $K$-theory.

We shall call the  generalized symmetric product functor with points having labels, or ``multiplicities'', in a partial monoid $M$, the Dold-Thom functor with coefficients in $M$.
Certainly, the construction of such functors via $\Gamma$-spaces is most appealing. However, if we start with a spectrum, constructing the corresponding Dold-Thom functor using $\Gamma$-spaces is not an entirely straightforward procedure since the $\Gamma$-space naturally associated to a spectrum is not injective. The purpose of the present note is to show how a connective spectrum $\mathcal M$ gives rise to an explicit partial monoid $M$ such that the homotopy of the Dold-Thom functor with coefficients in $M$ coincides, as a functor, with the homology with coefficients in $\mathcal M$. The construction is based on a trivial observation: if $Y$ is a space and $X$ is a pointed space, the space of maps from $Y$ to $X$ has a commutative partial multiplication with a unit.

\section{Partial monoids and infinite loop spaces}
\subsection{Partial monoids.}
Most of the following definitions appear in \cite{Segal:Inventiones}. A  partial monoid
$M$ is a topological space equipped with a subspace $M_{(2)}\subseteq M\times M$ and an addition map $M_{(2)}\to M$, written as $(m_1,m_2)\to m_1 + m_2$, and satisfying the following two conditions:
\begin{itemize}
\item there exists $0\in M$ such that $0+ m$ and $m+0$ are
defined for all $m\in M$ and such that $0+ m=m+0=m$;
\item for all $m_1,m_2$ and $m_3$ the sum $m_1+(m_2+ m_3)$
is defined whenever $(m_1+m_2)+ m_3$ is defined, and both are equal.
\end{itemize}
We shall say that a partial monoid is abelian  if for all $m_1$
and $m_2$ such that $m_1+ m_2$ is defined, $m_2+ m_1$ is also
defined, and both expressions are equal.

The
classifying space
$BM$ of a partial monoid $M$ is defined
as follows. Let $M_{(k)}$ be the subspace of $M^k$ consisting of
composable $k$-tuples. The $M_{(k)}$ form a simplicial space, with
the face operators $\partial_i: M_{(k)}\to M_{(k-1)}$ and the
degeneracy operators $s_i:M_{(k)}\to M_{(k+1)}$ defined as
$$
\begin{array}{rcll}
\partial_i(m_1,\ldots,m_k)&=&(m_2,\ldots,m_k)&\textrm{if\ } i=0\\
&=&(m_1,\ldots,m_i+ m_{i+1},\ldots,m_k)&\textrm{if\ } 0<i<k\\
&=&(m_1,\ldots,m_{k-1})&\textrm{if\ } i=k
\end{array}
$$
and
$$s_i(m_1,\ldots,m_k)=(m_1,\ldots,m_i,0,m_{i+1},\ldots,m_k)\quad
\textrm{if\ } 0\leq i\leq k.$$ The classifying space $BM$ is the
realization of this simplicial space. In the case when $M$ is a
monoid, $BM$ is its usual classifying space. If $M$ is a partial
monoid with a trivial multiplication (that is, the only composable
pairs of elements in $M$ are those containing $0$), the space $BM$
coincides with the reduced suspension $\Sigma M$.

A homomorphism $f: M \to N$ is a map such that whenever $m_1+m_2$ is defined,
$f(m_1)+f(m_2)$ is also defined and equal to $f(m_1+m_2)$.
If $f$, considered as a map of sets, is an inclusion, we say that $M$ is a 
partial
submonoid 
of $N$.

\subsection{Partial monoids and spectra.}

Given a pointed space $X$ we shall write $\Omega X$ for the space all maps $\R\to X$ supported (that is, attaining a value distinct from the base point of $X$) inside a compact subset of $\R$. If $\Omega_{\epsilon}X$  denotes the space of all maps $\R\to X$ supported  in $[-\epsilon, \epsilon]$, with the compact-open topology,  $\Omega X$ is given the weak topology of the union $\cup_{\epsilon}\Omega_{\epsilon}X$. The usual loop space can be identified with $\Omega_1 X$ and the inclusion $\Omega_1 X\hookrightarrow\Omega X$ is a homotopy equivalence.

If $X$ is an abelian partial monoid and $Y$ is a space, the space of all continuous maps $Y\to X$ is also an abelian partial monoid. Two maps $f,g$ are composable in this partial monoid if at each point of $Y$ their values are composable. In particular, if $X$ is a pointed topological space, then, considering $X$ as a monoid with the trivial multiplication, we see that $\Omega X$ is an abelian partial monoid; two maps in $\Omega  X$ are composable if their supports are disjoint.

The space $\Omega  X$ as defined here is much better behaved with respect to this partial multiplication than the usual loop space: while a generic element of the usual loop space is only composable with the base point, each element of $\Omega  X$ is composable with a big (in a sense that we need not make precise here) subset of $\Omega  X$.

The partial multiplication on $\Omega^n X$ for $n>1$ can be defined inductively for all the pairs of maps from $\R$ to the partial monoid  $\Omega^{n-1} X$ whose values are composable at each point.  This is, of course, the same as treating the points of  $\Omega^n X$ as maps $\R^n\to X$ and defining the composable pairs of maps as those with disjoint supports in $\R^n$.

Let now $\mathcal M$ be a connective spectrum. We construct  a partial abelian multiplication on its infinite loop space as follows. First, let us replace inductively $\mathcal M_0$ by a point and the spaces $\mathcal M_i$ for $i>0$ by the mapping cylinders of the structure maps $$\Sigma \mathcal
M_{i-1}\to\mathcal M_i, $$ obtaining a spectrum $\mathcal M'$. This, in particular, allows us to assume that all the structure maps are inclusions. As a consequence, we have inclusion maps $$\Omega^{i-1} \mathcal M'_{i-1}\to\Omega^i \mathcal M'_i,$$ which send composable $k$-tuples of elements to composable $k$-tuples for all $k$. The union of all the spaces 
$$\Mi{i}=\Omega^i \mathcal M'_i$$ is the infinite loop space for the spectrum $\mathcal M$ and it naturally has the structure of a partial abelian monoid.


\subsection{Dold-Thom functors.}

Let $M$ be an abelian partial moniod and $X$ a topological space with the base point $*$. We define the configuration space $M_n[X]$ of at most $n$ points in $X$ with labels in $M$ as follows. For $n>0$ let $W_n$ be the subspace of the symmetric product $SP^n(X\wedge M)$ consisting of points $\sum_{i=1}^{n}(x_i,m_i)$ such that the $m_i$ are
composable; $W_0$ is defined to be a point. The space $M_n[X]$ is the quotient of $W_n$ by the relations
$$(x,m_1)+(x,m_2)+\ldots=(x,m_1+m_2)+\ldots,$$
where the omitted terms on both sides are understood to coincide,  and
$$(x, 0)=(*, m)$$
for all $x$ and $m$. This quotient map commutes with the inclusions of $W_n$ into
$W_{n+1}$ coming, in turn, from the inclusions $SP^n(X\wedge M)\to
SP^{n+1}(X\wedge M)$, and, hence, $M_n[X]$ is a subspace of
$M_{n+1}[X]$.

The Dold-Thom functor of $X$ with coefficients in $M$ is the
space
$$M[X]=\bigcup_{n>0}M_n[X]. $$

The Dold-Thom functor with coefficients in the monoid of
non-negative integers is the infinite symmetric product. If $M$ has
trivial multiplication, we have $$M[X]=M_1[X]=X \wedge M.$$

The composability of labels in a configuration is essential for the
functoriality of $M[X]$. A based map $f:X\to Y$ induces a map
$M[f]:M[X]\to M[Y]$ as follows: a point $\sum(x_i,m_i)$ is sent to
the point $\sum(y_j,n_j)$ where the label $n_j$ is equal to the sum
of all the $m_i$ such that $f(x_i)=y_j$.

Apart from the infinite symmetric products, Dold-Thom functors
generalize classifying spaces: for any partial monoid $M$ its
classifying space $BM$ is homeomorphic to $M[S^1]$. To construct the
homeomorphism, take the lengths of the intervals between the
particles to be the barycentric coordinates in the simplex in $BM$
defined by the labels of the particles. Similarly, the classifying
space of an arbitrary $\Gamma$-space can be constructed in this way
(modulo some technical details), see Section~3 of
\cite{Segal:Topology}. The identification of $BM$ with $M[S^1]$ also
makes sense for non-abelian monoids. The construction of a
classifying space for a monoid as a space of particles on $S^1$ was
first described in \cite{McCord}.

In the  case when  $M$ is a partial abelian monoid coming from a spectrum $\mathcal M$ it is convenient to consider a different topology on $M_n[X]$. Namely, we define $M_n[X]$ as the the union of the spaces $\Mi{i}_n[X]$ with the weak topology. Then, as above, the Dold-Thom functor with coefficients in $M$ associates to a space $X$ the space $M[X]=\cup\, M_n[X]$. 
The technical advantage provided by this modification is that for any compact space $Y$ the map $Y\to M[X]$ factorizes through $\Mi{i}_n[X]$ for some $i$ and $n$.

\medskip

We have the following generalization of the Dold-Thom theorem:

\begin{theorem}\label{DT}
Let $M$ be partial abelian monoid corresponding to a connective spectrum $\mathcal M$. Then the spaces $M[S^n]$ form a connective spectrum weakly equivalent to $\mathcal M$. The functor
$$X\to\pi_*M[X]$$ coincides 
with the reduced homology with coefficients in $\mathcal M$.
\end{theorem}

\subsection{Other partial monoids}
The subject of this note are the partial monoids coming from spectra. Nevertheless, it is worth pointing out that there are many examples of partial abelian monoids, apart from infinite loop spaces, such that the corresponding Dold-Thom functors are linear. It is not hard to prove using the methods similar to those of Dold and Thom that if a partial monoid $M$ has a good base point and if the complement of the subspace of composable pairs has, in a certain sense, infinite codimension in $M\times M$, then the corresponding Dold-Thom functor defines a homology theory. We have already mentioned Segal's partial monoid of vector subspaces of $\R^{\infty}$ defined
in \cite{Segal:K-theory}. Among other examples we have:

\medskip

{\bf 1.} The configuration space of several (between 0 and $\infty$) distinct points in $\R^{\infty}$, with the sum of two disjoint configurations being defined as their union. For configurations with points in common the sum is not defined. The unit is the point $\emptyset$ thought of as the configuration space of 0 points. More generally, one can consider the configuration spaces of distinct particles in $\R^{\infty}$, labelled by points of a fixed space $M$. This was done in the paper \cite{Sh} by K.\ Shimakawa; the homology theory produced by this construction assigns to $X$ the stable homotopy of $X\wedge M$. (In fact, Shimakawa considers a more general situation of configurations in $\R^{\infty}$ with partially summable labels belonging to some partial abelian monoid. Such configuration space is then itself a partial abelian monoid which, as Shimakawa proves, always gives rise to a homology theory.)

\medskip

{\bf 2.} The space of all $n$-dimensional closed compact submanifolds of $\R^{\infty}$ with the sum being the union of the submanifolds, whenever they do not intersect. This construction, however, adds nothing substantially new to the previous example. Indeed, since the dimension of submanifolds is finite and their codimension is infinite, their connected components can be shrunk in size simultaneously (at least in a compact family of submanifolds) and we see that the partial monoid of $n$-submanifolds of $\R^{\infty}$ is weakly homotopy equivalent to the labelled configuration space of particles in $\R^{\infty}$, with labels in $\cup_{M} B{\rm Diff}(M)$, the space of all connected $n$-submanifolds of $\R^{\infty}$.

\medskip

{\bf 3.} The space of all (piecewise smooth) spheres in $\R^{\infty}$. The operation is the join of two spheres inside $\R^{\infty}$, and it is defined whenever any two intervals connecting points of the two summands are disjoint.


\section{Properties of Dold-Thom functors coming from spectra}

Given a partial abelian monoid $M$, a subset $Z\subset M$, and a
homotopy $$s_t:M\to M,$$ with $t\in [0,1],$ we say that $s_t$ is a
{\em deformation of $M$, admissible with respect to $Z$} if
\begin{itemize}
\item $s_0=Id$ and $s_t(0)=0$ for all $t$;
\item for any set of composable
elements $m_i\in M$, the set $s_t(m_i)$ is also composable for all
$t$;
\item if a set of composable
elements $m_i\in M$ is composable with $m'\in Z$, then the set
$s_t(m_i)$ is composable with $m'$ for all $t$.
\end{itemize}

\medskip

In what follows $M=\cup_i\, \Mi{i}$ is a partial abelian monoid  coming from a connective
spectrum.

\begin{lemma}\label{almost} 
For each $i$ and for any compact subset $Z\subset \Mi{i}$ there exists a deformation
$$d^i_t:\Mi{i}\times[0,1]\to \Mi{i+1},$$ admissible with respect to $Z$, such that any
element of $d^i_1(\Mi{i})$ is composable in $\Mi{i+1}$ with any element of
$Z$.
\end{lemma}

\begin{proof}
Consider the points of $Z$ and those of  $\Mi{i}$ as maps of $\R^{i+1}$ to $\mathcal M'_{i+1}$; let $x_1,\ldots,x_{i+1}$ be the coordinates in $\R^{i+1}$. Since $Z$ is compact there exists $a\in\R$ such that $f\in Z$ implies that the support of $f$ is contained in the half-space $x_{i+1}<a$.

Now, define $d^i_t$ for $0\leq t\leq 1$ by
$$d^i_t(f)(x_1,\ldots,x_{i+1})=
\begin{cases}
f\left(x_1,\ldots ,x_i, x_{i+1}-\frac{1}{x_{i+1}-a-1+t^{-1}}\right)\  {\text{if}} \ x_{i+1}>a+1-t^{-1};\\
\text{base point in}\  \mathcal M'_{i+1}\ \text{if}\ x_{i+1}\leq a+1-t^{-1}.
\end{cases}$$
The support of each element of $d^i_1(\Mi{i})$ is contained in
the half-space $x_{i+1}>a$, therefore each element of
$d^i_1(\Mi{i})$ is composable with each element of $Z$ in
$\Mi{i+1}$. Also, the deformation $d^i_t$, as a deformation of $\Mi{i}$,
is admissible with respect to $Z$. Indeed,  $d^i_t$ does
not change the support of $f$, and the composability of two elements
of $\Mi{i}$ only depends on whether their supports are disjoint or
not.
\end{proof}

\begin{lemma}\label{almost2} Let ${\bf I}$ be the abelian monoid whose elements are points of $[0,1]$ with the sum of two numbers being their maximum.
For each $i$ there is homomorphism of partial monoids $$h_i: \Mi{i}\to{\bf I}$$ which only vanishes at zero, and a deformation $$u_i: \Mi{i}\times [0,1]\to \Mi{i},$$ which  is constant on the set $h_i=1$,  decreases the value of $h_i$ strictly monotonically on the set $h_i<1$, at the value of the parameter $t<1$ is a homeomorphism and at $t=1$ retracts the subspace $h_i<1/4$ into the base point. Moreover $u_i$ and the restriction of $u_{i+1}$ to $\Mi{i}$ are homotopic.  
\end{lemma}

\begin{proof}
First we define inductively a collection of functions $l_i:\mathcal M'_i\to [0,1]$. The space $\mathcal M'_0$ is a point and we set $l_i$ to be equal to zero on it. Assume that we have already defined the function $l_k$. By construction, the space  $\mathcal M'_{k+1}$ is the mapping cylinder of the map $\Sigma \mathcal M'_k\to  \mathcal M_{k+1}$ induced by the structure map of the spectrum  $\mathcal M$. Consider the reduced suspension $\Sigma \mathcal M'_k$ as the product $ \mathcal M'_k\times [-1,1]$ with the identifications 
$$(*,s)\sim (x,1)\sim (x,-1)$$ 
for all $x\in \mathcal M'_k$ and $s \in [-1,1]$, and let $\tau$ be the cylinder coordinate, which is equal to 0 on  $\Sigma \mathcal M'_k$ and to 1 on $\mathcal M_{k+1}$. Then we set $l_{k+1}=1$ if $\tau=1$ and 
$$l_{k+1}((x,s), \tau)=\max\left( \min{(l_k(x), 2- 2|s|)}, \tau\right)$$
for $(x,s)\in \Sigma \mathcal M'_k$ and $\tau<1$. 

In the same vein, define a collection of retractions $w_i:\mathcal M'_i\times [0,1] \to \mathcal M'_i$. Let $q_t$ with $t\in [0,1]$ be a continuous family of continuous monotonic functions from $[0,1]$ to itself such that 
\begin{itemize}
\item $q_t(0)=0$ and $q_t(a)=a$ for all $t$ and $a\geq 1/2$;
\item $q_t$ is strictly monotonic for $t<1$, and $q_1(a)=0$ for $a<1/4$;
\item $q_t(a)>q_{t'}(a)$ for all $0<a<1/2$ and $t<t'$.    
\end{itemize}
Also define $r(\tau, l)$ as
$$r(\tau,l)=
\begin{cases}
\begin{array}{rll}
0&  {\text{if}} & 0\leq \tau \leq 1/2,\ 0\leq l \leq 1/2;\\
2\tau-1& {\text{if}} & 1/2 <  \tau \leq 1,\ 0\leq l \leq 1/2;\\
2\tau l - \tau & {\text{if}} & 0\leq \tau \leq 1/2,\ 1/2< l \leq 1;\\
-2\tau l +3\tau +2l -2 &  {\text{if}} & 1/2 <  \tau \leq 1,\ 1/2< l \leq 1.
\end{array}
\end{cases}$$

Then, assuming that we have already defined $w_k$, we set $w_{k+1}$ to be constant on $\mathcal M_{k+1}$ and on the set $\tau<1$ we define 
$$(w_{k+1})_t((x,s),\tau)=\bigl( ((w_k)_t(x), 1-q_t(1-s)), tr(\tau,l) + (1-t)\tau\bigr)$$
if $s\geq 0$, and
$$(w_{k+1})_t((x,s),\tau)=\bigl( ((w_k)_t(x), -1+q_t(1+s) ), tr(\tau,l) + (1-t) \tau \bigr)$$
for negative $s$. 

It is verified directly that $w_i$ is constant on the set $l_i=1$,  decreases the value of $l_i$ strictly monotonically on the set $l_i<1$, at the value of the parameter $t<1$ is a homeomorphism and at $t=1$ deforms the subspace $l_i<1/4$ into the base point. Now, define the function $h_i$ on $\Mi{i}$ on a map $\alpha :\R^i\to \mathcal M'_i$ as the maximal value of $h_i\circ \alpha$ and the deformation $u_i$ of  $\Mi{i}$ as that induced by the deformation $w_i$. It is clear that the conditions of the lemma are then satisfied.

\end{proof}

\begin{lemma}\label{almost3} The set $\pi_0 (M)$ is an abelian group with the addition induced by the partial addition in $M$.
\end{lemma}

\begin{proof}
It is sufficient to notice that if points of $M$ are thought of as maps of $\R$ to $\Omega^{i-1}\mathcal M'_i$ for some $i$, an inverse to a map $\gamma$ in $\pi_0 M$ will be given by $\gamma(C-t)$ for sufficiently big $C$.
\end{proof}

For a space $X$ let $M_{<a>}[X]\subset M[X]$ be the subset of configurations whose
coefficients are composable with $a\in M$.

\begin{lemma}\label{plus} The inclusion of
$M_{<a>}[X]$ into $M[X]$ induces isomorphisms on all homotopy
groups. Moreover, if $\sum (x_i, m_i) \in M_{<a>}[X]$, the map
$$M_{<a+\sum m_i>}[X]\to M_{<a>}[X]$$
given by summing with $\sum (x_i, m_i)$, is also a weak homotopy equivalence.
\end{lemma}

\begin{proof}
We first need to show that the image of any map $f$ of a finite cell complex $Y$ into $M[X]$ can be deformed into $M_{<a>}[X]$. Assume that the image of $f$ is contained in $\Mi{i}[X]$. Then, applying the deformation $d^i_t$ from Lemma~\ref{almost}, with $Z=\{a\}$, to all the labels of the configurations $f(y)$, where $y\in Y$, we get a homotopy of $f$ to a map of $Y$ into $M_{<a>}[X]$.

In order to establish the second claim, for each $m_i$ choose $\tilde{m}_i$ in such a way  that the class of $\tilde{m}_i$ is inverse to that of $m_i$ in $\pi_0 M$ and so that the sum $a+\sum m_i + \sum \tilde{m}_i$  is defined.

It is then easy to see that adding $\sum ( x_i, m_i)$ and then $\sum ( x_i,  \tilde{m}_i)$ to a configuration in $M_{<a + \sum m_i +\sum \tilde{m}_i>}[X]$ is homotopic to the natural inclusion
$$M_{<a + \sum m_i + \sum \tilde{m}_i>}[X]\hookrightarrow M_{<a>}[X].$$
Since all the natural inclusions between the spaces $M_{<a>}[X]$ for different $a$ are homotopy equivalences, it follows that the map in the statement of the lemma is surjective on the homotopy groups. Replacing in this argument $a$ by $a+\sum m_i$ we see that the map
$$M_{<a + \sum m_i + \sum \tilde{m}_i>}[X]\hookrightarrow M_{<a+ \sum m_i >}[X],$$
given by adding $\sum ( x_i,  \tilde{m}_i)$, is also surjective on homotopy and this proves the lemma.
\end{proof}


\section{Quasifibrations of Dold-Thom functors}
The proof of Theorem~\ref{DT} is based on the original argument of Dold and Thom
\cite{DT}, see also \cite{AGP, Hat}. We have the following fact:
\begin{prop}\label{prop:DT}
Let $M$ be a 
partial monoid coming from a connective spectrum.
If $X$ is a cell complex and $A\subset X$ is a subcomplex, the map $M[X]\to M[X/A]$ is a quasifibration with the fibre $M[A]$.
\end{prop}

The rest of this section is dedicated to the proof of this statement. Note that we do not require $A$ to be connected.

\medskip

We shall need the following criterion for quasifibrations. Let $p: E\to B$ be a map which is  quasifibration over $B'\subset B$. Assume that for any  
compact $C\subset B$ there is a homotopy $f_t$ of the inclusion map $i: C\hookrightarrow B$ to a map $C\to B'$, which maps $C\cap B'$ to $B'$ for all $t$. Further, suppose that for any compact $\widetilde{C}\subseteq p^{-1}(C)$ there is a homotopy $\tilde{f}_t$ of the inclusion map $\widetilde{C}\to E$ to a map $\widetilde{C}\to p^{-1}(B')$ which maps $\widetilde{C}\cap p^{-1}(B')$ to $p^{-1}(B')$ for all $t$,  and such that $$p\circ \tilde{f}_t = f_t\circ p.$$
Moreover, assume that $f_t$ and $\tilde{f}_t$ are well-defined up to homotopy.

Take a point $b\in B$ and $\tilde{b}\in p^{-1}(b)$. Then we have two paths: from $b$ to $b'\in B'$ and from $\tilde{b}\in p^{-1}(b)$ to $\tilde{b}'\in p^{-1}(b')$, the latter covering the former. This gives a  map of the homotopy groups
$$\pi_i(p^{-1}(b),\tilde{b}) \to \pi_i(p^{-1}(b'), \tilde{b}').$$

\begin{lemma}\label{criterion}
Assume that all the above maps of homotopy groups of the fibres are isomorphisms for all $i\geq 0$. Then the map $p$ is a quasifibration.
\end{lemma}
This lemma is a version of Hilfssatz~2.10 of \cite{DT} and the proof is, essentially, the same.

\medskip

We denote by $p$ the projection map $$X\to X/A$$ and by $\pi$ the
induced map $$M[X]\to M[X/A].$$ We shall prove by induction
on $n$ that $\pi$ is a quasifibration over $M_{n}[X/A]$. This, by
Satz 2.15 of \cite{DT} (or Theorem A.1.17 of \cite{AGP}) will imply
that $\pi$ is a quasifibration over the whole $M[X/A]$.

Assume that $\pi$ is a quasifibration over $M_{n-1}[X/A]$. According
to Satz 2.2 of \cite{DT} (or Theorem A.1.2 of \cite{AGP}) it is
sufficient to prove that $\pi$ is a quasifibration over
$M_n[X/A]-M_{n-1}[X/A]$, over a neighbourhood of $M_{n-1}[X/A]$ in
$M_{n}[X/A]$ and over the intersection of this neighbourhood with
$M_n[X/A]-M_{n-1}[X/A]$.

It will be convenient to speak of  delayed homotopies. A
delayed homotopy is a map $f_t:A\times[0,1]\to B$  such that for
some $\varepsilon>0$ we have $f_t=f_0$ when $t\leq\varepsilon$. A
map $p:E\to B$ is said to have the delayed homotopy lifting property
if it has the homotopy lifting property with respect to all delayed
homotopies of finite cell complexes into $B$. It is clear that a map
that has the delayed homotopy lifting property is a quasifibration.

\begin{lemma}
Let $B$  be an arbitrary subspace of $M_n[X/A]-M_{n-1}[X/A]$. The map
$\pi$, when restricted to $\pi^{-1}(B)$, has the delayed homotopy
lifting property.
\end{lemma}
\begin{proof}

Let $$f_t:Z\times [0,1]\to B$$ be a delayed homotopy of a finite
cell complex $Z$ into $B$, such that $f_t=f_0$ for $t\leq\varepsilon$,
and let $$\tilde{f}_0:Z\to \pi^{-1}(B)$$ be its lifting at $t=0$.

Notice that $M_n[X/A]-M_{n-1}[X/A]$ can be thought of as the
subspace of $M_n[X]$ consisting of configurations of $n$ distinct
points, all outside $A$ and with non-trivial labels. Therefore, we
can think of $B$ as of a subspace of $M_n[X]$.

Define $g:Z\to M[X]$ as the difference $$g(z)=\tilde{f}_0(z)-f_0(z).$$
The map $g$ is well-defined, continuous and its image belongs to
$M[A]$. Since $Z$ is compact, the image of $\tilde{f}_0$ belongs to
$\Mi{i}[X]$ for some $i$; it follows that the coefficients of $g(z)$
are composable with the coefficients of $f_0(z)$ in $\Mi{i}$ for all
$z\in Z$.

Lemma~\ref{almost} guarantees the existence of a deformation $d^i_t$ of $\Mi{i}$ inside $\Mi{i+1}$ such that each point in $d^i_1(\Mi{i})$
is composable in $\Mi{i+1}$ with each point in the image of $Z\times
[0,1]$ under $f_t$. There is an induced deformation of $\Mi{i}[A]$
which we also denote by $d^i_t$.

Define the homotopy $g_t:Z\to M[A]$ as $d^k_{t\varepsilon^{-1}}\circ
g$ for $0\leq t<\varepsilon$ and as $d^k_1\circ g$ for
$\varepsilon\leq t\leq 1$.
Lemma~\ref{almost}  implies that $g_t$ is well-defined and is composable with
$f_t$ for all $t$. Consider the map $\tilde{f}_t=f_t+g_t:Z\to M[X]$.
By construction, it lifts $f_t$.
\end{proof}

It remains to see that the map $\pi$ is a quasifibration over some neighbourhood  of $M_{n-1}[X/A]$.

Recall from Lemma~\ref{almost2} that  there is a homomorphism of partial monoids
$$\Mi{i}\to {\bf I}$$
which sends a map $\alpha$ to the maximal value of $h_i\circ \alpha$. It gives rise to a map
$$v_i: \Mi{i}_n[X] \to {\bf I}_n[X].$$

If $X$ and $A$ are cell complexes, let  $\Delta$ be the subspace of ${\bf I}_n [X]$ consisting of configurations which either contain a point of $A$, or have less than $n$ points.  It is not hard to show that ${\bf I}_n[X]$ is a cell complex and $\Delta$ is a subcomplex. In particular, $\Delta$ is a strong deformation retract of its neighbourhood $U\subset {\bf I}_n[X]$. Let $$d_t: {\bf I}_n[X]\times [0,1]\to {\bf I}_n[X]$$ be the deformation that retracts $U$ to $\Delta$. We can assume that $d_t$ is a homeomorphism for all $t<1$. Then $d_t$ can be lifted to a deformation
$$\widetilde{D}_t: \Mi{i}_n[X] \times [0,1]\to \Mi{i}_n[X],$$
that retracts the open subset $v_i^{-1}(U)$ to the subspace $v_i^{-1}(\Delta)\subset \Mi{i}_n[X] $, which consists of configurations which either contain a point of $A$, or have less than $n$ points.

Now, by construction, there exists a deformation $${D}_t: \Mi{i}_n[X/A] \times [0,1]\to\Mi{i}_n[X/A],$$
such that $$D_{t}\circ \pi= \pi\circ\widetilde{D}_t$$ for all $t$. In particular, $D_t$ retracts the open neighbourhood $\pi(v_i^{-1}(U))$ of $\Mi{i}_{n-1}[X/A]$ onto $\Mi{i}_{n-1}[X/A]$.

Let $V$ be the union of all the neighbourhoods $\pi(v_i^{-1}(U))$ for all $i$ in $M[X/A]$. Then $V$ is open in $M[X/A]$ and it follows from Lemma~\ref{plus} that the projection $\pi^{-1}(V)\to V$ satisfies the conditions of Lemma~\ref{criterion}.


\section{On the spectrum $M[\mathcal S]$}
\subsection{The spectrum $M[\mathcal S]$ and the proof of Theorem~\ref{DT}.}
Proposition~\ref{prop:DT} with $X=D^n$ and $A=\partial D^n$ gives
the quasifibration
$$M[S^{n-1}]\to M[D^n]\to M[S^n].$$
The space $M[D^n]$ is contractible, and therefore, we have weak
homotopy equivalences $M[S^{n-1}]\simeq\Omega M[S^n]$ and the spaces
$M[S^i]$ for $i\geq 0$ form an $\Omega$-spectrum which we denote by
$M[\mathcal S]$.

More generally, given $X$, the cofibration $X\to
CX\to \Sigma X$ gives rise to a weak homotopy equivalence
$M[X]\simeq \Omega M[\Sigma X]$,
and given an inclusion map $i:A\hookrightarrow X$,
the cofibration $A\to Cyl(i)\to X\cup_i CA$ gives rise to
an exact sequence $$\ldots\pi_* M[A]\to \pi_*M[X] \to \pi_* M[X\cup_i CA]\to\pi_{*-1} M[A]\to\ldots.$$
Here $CX$ is the cone on $X$ and
$Cyl(i)$ is the cylinder of the map $i$. Since $\pi_* M[X]$ is,
clearly, a homotopy functor, this means that the groups $\pi_* M[X]$
form a reduced homology theory.

There is a natural transformation of the homology with coefficients
in $M[S^i]$ to $\pi_*M[X]$, induced by the obvious map
$M[S^{i}]\wedge X\to M[S^{i}\wedge X]$ that sends $\bigl(\sum
m_iz_i,x\bigr)$ to $\sum m_i(z_i,x)$:
$$\lim_{i\to\infty}\pi_{k+i}\bigl(M[S^{i}]\wedge X\bigr)\to
\lim_{i\to\infty}\pi_{k+i} M[S^{i}\wedge X]=\pi_k M[X].$$ If $X$ is
a sphere, the Freudenthal Theorem implies that this is an
isomorphism. Hence, $\pi_*M[X]$ coincides as a functor, on connected
cell complexes, with the homology with coefficients in $M[\mathcal
S]$.

\subsection{The weak equivalence of $M[\mathcal S]$ and $\mathcal M$.}

We shall first construct a weak homotopy equivalence between the
infinite loop spaces of the spectra  $\mathcal M$ and $M[\mathcal
S]$ and then show that there exists an inverse to this equivalence,
which is induced by a map of spectra. This will establish
Theorem~\ref{DT}.

Let $I$ be the interval $[-1/2,1/2]$. By Proposition~\ref{prop:DT}
the map $I\to S^1$ which identifies the endpoints of $I$ induces a
quasifibration $M[I]\to M[S^1]$ with the fibre $M\simeq \Omega^{\infty}\mathcal M_{\infty}$.

Since $M[I]$ is contractible, it follows that $M$ is
weakly homotopy equivalent to $\Omega M[S^1]$; the weak equivalence
is realized by the map $\psi$ that sends $m\in M$ to the
loop (parameterized by $I$) whose value at the time $t\in I$ is the
configuration consisting of one point with coordinate $-t$ and label
$m$.

Let us now define the map of spectra $\Phi$, inverse to the above weak equivalence $\psi$.

For $n>0$ identify $S^n$ with the $n$-dimensional cube $I^n=[-1/2,1/2]^n\subset \R^n$,
modulo its boundary. Fix a homeomorphism of the interior of $I^n$ with $\R^n$,
say, by sending each coordinate $u_k$ to $\tan{\pi u_k}$. Think of a point in $M[S^n]$
as of a sum $\sum (x_\alpha, m_\alpha)$ where $x_\alpha\in I^n$ and $m_\alpha$ are maps from $I^n$ to $\mathcal M'_n$. Assume that the maps $m_\alpha(-x_\alpha)$ are composable as maps from some $\R^i$ to $\mathcal M'_{n+i}$. Then their sum is a well-defined point of $\mathcal M'_n$ which does not depend on $i$. We set $\Phi_n(\sum (x_\alpha, m_\alpha))$ to be equal to this sum.

The map $\Phi_n: M[S^n]\to \mathcal M'_n$
is only partially defined, but this problem can be circumvented as follows.
Let $M^*[S^n]\subset M[S^n]$ be the subspace consisting of the configurations $\sum (x_\alpha, m_\alpha)$ such that for some $i$ the points $m_\alpha(q_\alpha)\in\Omega'^i\mathcal M'_{n+i}$ are composable (that is, have disjoint supports in $\R^i$) for any choice of the $q_{\alpha}\in I^n$.

The spaces $M^*[S^n]$ form a sub-spectrum $M^*[\mathcal S]$ of
$M[\mathcal S]$. Indeed, the structure map of $M[\mathcal S]$ sends
$\sum (x_{\alpha}, m_{\alpha})\in M[S^n]$ to the loop $t\to\sum
((x_{\alpha},t), m_{\alpha})$, where $t\in I$. The map $\Phi_n$ is well-defined on $M^*[S^n]$.
If we define the map $\Phi_0$ simply as $\Omega\Phi_1$,  it is clear
from the construction that $\Phi_0\circ\psi$ is the identity map on
$\mathcal M_0$.  The proof will be finished as soon as we prove the following
\begin{prop}\label{prop:eqspec}
The inclusion map $M^*[S^n]\to M[S^n]$ is a weak homotopy equivalence for all $n>0$.
\end{prop}
\begin{proof}
Define $M_k^*[S^n]$ as $M^*[S^n] \cap M_k[S^n]$. It sufficient to prove that the inclusion $M_k^*[S^n]\to M_k[S^n]$ is a weak homotopy equivalence for all $k$.

For $k=1$ this inclusion is the identity map. Also, for $k>1$ the map 
\begin{equation}\label{smth}
M_k^*[S^n]/ M_{k-1}^*[S^n]\to M_k[S^n]/M_{k-1}[S^n]
\end{equation}
is a weak homotopy equivalence. 
Indeed, take a map $f: S^j\to M_k[S^n]/ M_{k-1}[S^n]$ and let us show that the labels of the points in each configuration in the image of $f$ can be pushed away from each other, thus deforming the image of $f$ into $M_k^*[S^n]/ M_{k-1}^*[S^n]$.

By forgetting the labels in the configurations we get a map 
$$(S^j-f^{-1}(*) )\to B_k(S^n)$$ 
to the configuration space of $k$ distinct particles in $S^n$. Without loss of generality we can assume that the boundary of $S^j-f^{-1}(*)$ is a collared (for instance, smooth) hypersurface in $S^{j}$, so that this map gives rise to a bundle $\xi$ of $k$-element sets over the compactification $C$ of $S^j-f^{-1}(*)$. In turn, this bundle of sets gives rise to a $k$-vector bundle $\eta$ whose fibre is spanned by the elements of the corresponding fibre of $\xi$. Since $C$ is compact, $\eta$ can be considered as a subbundle of some trivial bundle. What this means is that there exists an $N$ such that given a configuration in the image of $f$ we can assign a unit vector in some $\R^N$ to each of its points so that the vectors assigned to all the points are mutually orthogonal.

Now, the image of $f$ is contained in $\Mi{i}_k[S^n]/\Mi{i}_{k-1}[S^n]$ for some $i$. We can treat the labels of the configurations in the image of $f$ as elements of 
$$\Mi{i+N}=\Omega^i\Omega^N\mathcal M'_{i+N},$$ and write them as functions $m(x,y)$ with $x\in\R^i$ and $y\in\R^N$. For $z\in C$ write $f(z)=\sum (q_{\alpha}, m_\alpha)$ where $q_\alpha\in S^n$ and $m_\alpha\in \Mi{i+N}$, and define
$f_t(z)$ as  $$f_t(z)=\sum (q_{\alpha}, m_\alpha(x, y+tv_{\alpha})),$$
where $v_\alpha$ are the orthogonal vectors associated to the points $q_{\alpha}$ of the configuration $f(z)$. 
Then for $t$ sufficiently large, the image of $f_t$ will lie in $M_k^*[S^n]/ M_{k-1}^*[S^n]$.

The subspaces $M_{k-1}[S^n]\subset M_k[S^n]$ and $M_{k-1}^*[S^n]\subset M_k^*[S^n]$ are not, strictly speaking, neighbourhood deformation retracts, but each of these subspaces has a neighbourhood such that the map of any compact into this neighbourhood can be retracted onto the subspace. This is sufficient to claim that the fact that the maps (\ref{smth}) are weak homotopy equivalences implies that the inclusion $M_k^*[S^n]\to M_k[S^n]$ induce isomorphisms in homology. The fundamental groups of these two spaces are easily seen to be abelian for $k>1$, just as in the case of the usual symmetric products, and, hence, this inclusion is a weak homotopy equivalence.

\end{proof}

\subsection*{Acknowledgments}
I would like to thank Ernesto Lupercio for conversations that
provided the motivation for this paper, and to Norio Iwase, Dai
Tamaki and Peter Teichner for pointing out useful references. I am
also grateful to Max-Planck-Institut f\"ur Mathematik, Bonn for
hospitality during the stay at which this paper was written. This
work was partially supported by the CONACyT grant no.\ 44100.

\end{document}